\documentclass[10pt]{article}
\scrollmode

\usepackage{amsfonts}
\usepackage{amsmath}
\usepackage{amssymb}
\usepackage{graphicx}

\allowdisplaybreaks

\def\eps{\varepsilon}
\def\qed{\hfill\rule{.2cm}{.2cm}}

\def\E{{\mathbb E}}
\def\Z{{\mathbb Z}}

\def\R{{\mathbb R}}

\def\1{{\mathbf 1}}

\def\a{\alpha}
\def\o{\omega}

 \def\O{\Omega}

\newcommand{\CC}          {\mathcal{C}} 
\newcommand{\BB}          {\mathcal{B}}

\newcommand{\LL}         {\mathcal{L}}

\newtheorem{theo}{Theorem}[section]

\newtheorem{lm}[theo]{Lemma}

\newtheorem{rmk}[theo]{Remark}

\def\beq{\begin{equation}}
\def\eeq{\end{equation}}
\newcommand{\bei}{\begin{itemize}}
\newcommand{\eei}{\end{itemize}}
\newcommand{\ben}{\begin{enumerate}}
\newcommand{\een}{\end{enumerate}}
\newcommand{\beqn}{\begin{eqnarray}}
\newcommand{\beqnn}{\begin{eqnarray*}}
\newcommand{\eeqn}{\end{eqnarray}}
\newcommand{\eeqnn}{\end{eqnarray*}}
\newcommand{\brm}{\begin{rmk}}
\newcommand{\erm}{\end{rmk}}

\title{Quenched invariance principles for random walks on percolation clusters. } 
\author{ 
P.~Mathieu \footnote{ Universit\'e de Provence, CMI, 39 rue Joliot-Curie, 13013 Marseille, FRANCE.  
pierre.mathieu@cmi.univ-mrs.fr}
\and 
A.~Piatnitski \footnote{Lebedev Physical Institute of Russian Academy of Sciences and 
Narvik Institute of Technology, P. O. Box 385, N-8505 Narvik, NORWAY. 
andrey@sci.lebedev.ru}      
}

\begin{document}
\maketitle

\begin{abstract} 
We prove the almost sure ('quenched') invariance principle for a random walker 
on an infinite percolation cluster in $\Z^d$, $d\geq 2$. 
\end{abstract}

%\vskip 3mm
%\noindent{\bf Keywords:} 
%\vskip 3mm
%\noindent{\bf AMS Classification numbers: } Primary: 

%%%%%%%%%%%%%%%%%%%%%%%%%%%%%%%%%%%%%%%%%%%%%%%%%%%%%%%%%%%%%%%%%
%%%%%%%%%%%%%%  Intro  %%%%%%%%%%%%%%%%%%%%%%%%%%%%%%%%%%%%%%%%%%
%%%%%%%%%%%%%%%%%%%%%%%%%%%%%%%%%%%%%%%%%%%%%%%%%%%%%%%%%%%%%%%%%

\section{Introduction} 

Consider super critical Bernoulli bond percolation in $\Z^d$, $d\geq 2$: 
for $x, y \in \Z^d$, we write: $x \sim y$ if $x$ and $y$ are neighbors in 
the grid $\Z^d$, and let 
$\E_d$ be the set of non-oriented nearest pairs $(x,y)$. 
We identify a sub-graph of $\Z^d$ with an  application $\omega\in\{0,1\}^{\E_d}$, 
writing $\omega(x,y)=1$ if the edge $(x,y)$ is present in $\omega$ and 
$\omega(x,y)=0$ otherwise. Thus $\Omega=\{0,1\}^{\E_d}$ is 
the set of sub-graphs of $\Z^d$. 
Edges pertaining to $\o$ are then called {\it open}. 
Connected components of such a sub-graph will 
be called {\it clusters} and the cluster of $\omega$ containing a point 
$x\in\Z^d$ is denoted with $\CC_x(\o)$.

Define now $Q$ to be the probability measure on $\{0,1\}^{\E_d}$ under which  
the random 
variables $(\omega(e),\,e \in \E_d)$ are Bernouilli$(p)$ independent variables  
and let 
\beqnn p_c = \sup \{p\,;\,Q[\# \CC_0(\o)=\infty]= 0 \}
\eeqnn 
be the critical probability. 
It is known that $p_c\in]0,1[$, see \cite{kn:G}. 
Throughout the paper, we choose a parameter $p$ such  that 
\beq \label{p>p_c} p>p_c. 
\eeq 
Then, $Q$ almost surely, the graph $\o$ has a unique infinite cluster there after denoted with $\CC(\o)$. 

We are interested in the behaviour of the simple symmetric random walk on 
$\CC_0(\o)$: let $D(\R_+,\Z^d)$ be the space of c\`ad-l\`ag $\Z^d$-valued functions 
on $\R_+$ and $X(t)$, $t\in\R_+$, be the coordinate maps from $D(\R_+,\Z^d)$ to 
$\Z^d$. $D(\R_+,\Z^d)$ is endowed with the Skorohod topology. 
For a given sub-graph $\omega\in\{0,1\}^{\E_d}$, and for $x\in\Z^d$, let 
$P^\o_x$ be the probability measure on $D(\R_+,\Z^d)$ under which the coordinate process 
is the Markov chain starting at $X(0)=x$ and with generator 
\beqnn 
\LL^\o f(x)=\frac 1{n^\o(x)}\sum_{y\sim x} \o(x,y) (f(y)-f(x))\,, 
\eeqnn 
where $n^\o(x)$ is the number of neighbors of $x$ in the cluster $\CC_x(\o)$. 

The behaviour of $X(t)$ under $P^\o_x$ can be described as follows: 
starting from point $x$, the random walker waits for an exponential time of parameter $1$ and 
then chooses, uniformly at random,  
one of its neighbors in $\CC_x(\o)$, say $y$ and moves to $y$. 
This procedure is then iterated with independent hoping times. 
The walker clearly never leaves the cluster of $\o$ it started from. 
Since edges are not oriented, the  measures with weights $n^\o(x)$ on the 
possibly different clusters of $\o$ are reversible. 

Let $Q_0$ be the conditional measure $Q_0(.)=Q(.\vert \#\CC_0(\o)=\infty)$ and let 
$Q_0.P_x^\o$  be  the so-called {\it annealed} semi-direct product measure   
law defined by 
\beqnn Q_0.P_x^\o[\,F(\o,X(.))\,]=\int P_x^\o[\,F(\o,X(.))\,]\,dQ_0(\o)\,.\eeqnn  
Note that $X(t)$ is not Markovian  anymore under $Q_0.P_x^\o$. 
From \cite{kn:DFGW}, it is known that, under $Q_0.P_0^\o$, the 
process $(X^\eps(t)=\eps X(\frac t{\eps^2}),t\in\R_+)$ satisfies an invariance 
principle as $\eps$ tends to $0$ i.e. it converges in law to the law of a non-degenerate 
Brownian motion. The proof is based on the {\it point of view of the 
particule}. It relies on the fact that the law of the environment $\o$, 
viewed from the current position of the Markov chain is reversible, when considered 
under the annealed measure. It does not give any information on the behaviour 
of the walk for a typical choice of $\o$. 
On the other hand, only partial results in dimension higher than $4$ 
have been obtained for almost sure, 
also called {\it quenched}, invariance principles in the joint work of V.~Sidoravicius 
and A-S.~Sznitman, \cite{kn:SS}. Our result holds for any dimension:  

\begin{theo}
\label{theo}
$Q$ almost surely on the event $\#\CC_0(\o)=\infty$, under $P^\o_0$, the process 
$(X^\eps(t)=\eps X( t/\eps^2),t\in\R_+)$ converges in law 
as $\eps$ tends to $0$ to a 
Brownian motion with covariance matrix $\sigma^2Id$ where $\sigma^2$ 
is positive and 
does not depend on $\o$.  
\end{theo}

Our strategy of proof follows the classical pattern introduced by 
S.M.~Kozlov for averaging random walks with random conductances 
in \cite{kn:Ko}. 
The method of Kozlov was successfully used under ellipticity assumptions 
that are clearly not satisfied here. We refer in particular to  
first part of \cite{kn:SS} where random walks 
in elliptic environments are considered.  
The main idea is to modify the process $X(t)$ by the addition 
of a {\it corrector} in such a way that the sum is a martingale under 
$P^\o_0$ and use a martingale invariance principle. 
Then one has to prove that, in the rescaled limit, the corrector can be  
neglected or equivalently that the corrector has sub-linear growth. 
For this second step, 
in a classical elliptic set-up, one would invoke the Poincar\'e inequality 
and the compact embedding of $H^1$ into $L^2$. For percolation models, 
a weaker, but still suitable form of the Poincar\'e inequality was proved 
in the paper of 
P.~Mathieu and E.~Remy \cite{kn:MR}, see also \cite{kn:Ba}. 
However another difficulty arises: our reference measure is the counting 
measure on the cluster at the origin. When rescaled, it does converge to Lebesgue measure 
on $\R^d$ but, for a fixed $\eps$ it is of course singular. 
Thus rather than using classical functional analysis tools, one has to turn to 
$L^2$ techniques in varying spaces  or $2$ scale convergence arguments 
as they have been recently developped for 
the theory of homogenization of singular random structures in the work of A.~Piatnitski and 
V.~Jikov, see \cite{kn:JP}. 
An elementary self-contained construction of the corrector is given in section 
\ref{ssec:constr}. We also provide an approach to $2$ scale convergence avoiding explicit 
reference to the results of \cite{kn:JP}. 
For background material on homogenization theory in both periodic and random 
environments we refer to the book \cite{kn:JKO} where 
percolation models are considered in chapter 9. 

{\it Note on the constants: } throughout the paper $\beta$ and $c$ will denote  positive constants  
depending only on $d$ and $p$ whose values might change from place to place.

{\it Note: } N. Berger and M. Biskup recently obtained a proof of Theorem (\ref{theo}), 
see \cite{kn:BB}. Although they also rely on the construction of a corrector, 
their method  to prove the sub-linear growth  of the corrector is quite different from ours.

\vfill\eject 
%%%%%%%%%%%%%%%%%%%%%%%%%%%%%%%%%%%%%%%%%%%%%%%%%%%%%%%%%%%%%%%%%%
%%%%%%%%%%%%%% PART II %%%%%%%%%%%%%%%%%%%%%%%%%%%%%%%%%%%%%%%%%%%
%%%%%%%%%%%%%%%%%%%%%%%%%%%%%%%%%%%%%%%%%%%%%%%%%%%%%%%%%%%%%%%%%%

\vskip 1cm
\section{Proof of the theorem}
\setcounter{equation}{0}
\label{sec:proof}

%%%%%%%%%%%%%%%%%%%%%%%%%%%%%%%%%%%%%%%%%%%%%%%%%%%%%%%%%%%%%%%%%
%%%%%%%%%%%%%%  Intro  %%%%%%%%%%%%%%%%%%%%%%%%%%%%%%%%%%%%%%%%%%
%%%%%%%%%%%%%%%%%%%%%%%%%%%%%%%%%%%%%%%%%%%%%%%%%%%%%%%%%%%%%%%%%

Let $\vert x\vert=\max \vert x_i\vert$. 
We use the notation $x\cdot y$ for the scalar product of the two vectors $x,y\in\R^d$. 
We also use the notation $Q_0(.)=Q(.\vert\#\CC_0(\o)=\infty)$. 

\subsection{Tightness} \label{ssec:tightness} 

We start recalling the Gaussian upper bound obtained by M.~Barlow for walks 
on percolation clusters, see \cite{kn:Ba}. The corresponding lower 
bound also holds, but we won't need it here. Note that, Barlow's bound 
is only used in the proof of the tightness. Remember that $p>p_c$ so 
that, $Q$.a.s. the percolation sub-graph $\o$ contains a unique infinite cluster 
denoted with $\CC(\o)$. 

Statement from \cite{kn:Ba}: $Q$.a.s., for any $x\in\CC(\o)$ there exists a 
random variable $S_x$ such that whenever $x$ and $y$ belong to $\CC(\o)$ 
and $t\geq S_x$ then 
\beq \label{eq:Barlow} 
P^\o_x[X(t)=y]\leq ct^{-d/2} \exp(-\frac{\vert y-x\vert^2}{ct})\,. 
\eeq
Moreover, 
\beq\label{eq:Barlow2} 
Q[x\in\CC(\o)\,,S_x\geq t]\leq c\exp(-c t^{\,\eps(d)}) \hbox{ \, with\, $\eps(d)>0$}\,. 
\eeq 

((\ref{eq:Barlow}) is only stated with the further restriction that $t\geq\vert x-y\vert$ in 
\cite{kn:Ba}. If $t\leq\vert x-y\vert$, then (\ref{eq:Barlow}) follows from the Carne-Varopoulos 
bound, see \cite{kn:MR}, appendix C.)

\begin{lm} 
\label{tightness} 
$Q$ almost surely on the event $\#\CC_0(\o)=\infty$, under $P^\o_0$, the 
sequence of processes  
$(X^\eps(t)=\eps X(\frac t{\eps^2}),t\in\R_+)$ is tight in the Skorohod topology. 
\end{lm} 

{\it Proof}: 
it is sufficient to check that $Q$.a.s. on the event $\#\CC_0(\o)=\infty$,  for any $T>0$ one has 
\beqnn 
\limsup_{\delta\rightarrow 0}\, 
\limsup_{\eps\rightarrow 0}\, 
\sup_{\tau} \, E^\o_0[\vert X^\eps(\tau+\delta)-X^\eps(\tau)\vert^2]=0\,,
\eeqnn  
where $\tau$ is any stopping time in the filtration generated by $X^\eps$ that is bounded by $T$. 
See \cite{kn:EK}, page 138. 

But, using the Markov property, we get that for large enough $K$,  
\beqnn 
&&E^\o_0[\vert X^\eps(\tau+\delta)-X^\eps(\tau)\vert^2]\\ 
&=& \eps^2 E^\o_0[\vert X(\frac{\tau+\delta}{\eps^2})-X(\frac{\tau}{\eps^2})\vert^2]\\ 
&\leq& \eps^2 \sup_{y\in\CC_0(\o)\,;\, \vert y\vert\leq  K/\eps^2} E^\o_y[\vert X(\frac{\delta}{\eps^2})\vert^2] 
+e^{-K/2}\,,
\eeqnn 
In the 'sup', the restriction  $\vert y\vert\leq  K/\eps^2$ is justified 
since the walker makes more than $k$ steps in time $t$ with probability lower than $e^{-k/t}$. 
Since we are conditioning on the event $\CC_0(\o)=\CC(\o)$, one may replace the condition 
$y\in\CC_0(\o)$ by the condition $y\in\CC(\o)$ in the last term of this inequality. 

From (\ref{eq:Barlow}), it follows that the first term is bounded by $c\delta$ provided 
that $\frac \delta{\eps^2}\geq \sup_{y\in\CC(\o)\,;\, \vert y\vert\leq  K/\eps^2} S_y$. 
From (\ref{eq:Barlow2}), we get that 
\beqnn 
Q[\sup_{y\in\CC(\o)\,;\, \vert y\vert\leq  K/\eps^2} S_y>\frac \delta{\eps^2}]
\leq c \frac{ K}{\eps^2} \exp(-c (\frac \delta{\eps^2})^{\eps(d)})\,. \eeqnn 
Using the Borel-Cantelli lemma, we deduce that $Q$.a.s. on the event $\CC_0(\o)=\CC(\o)$ 
\beqnn \limsup_{\eps\rightarrow 0} 
\eps^2 \sup_{y\in\CC(\o)\,;\, \vert y\vert\leq  K/\eps^2} S_y=0\,, \eeqnn 
and the proof is completed by letting $K$ tend to $\infty$. \qed

\vskip1cm 
%%%%%%%%%%%%%%%%%%%%%%%%%%%%%%%%%%%%%%%%%%%%%%%%%%%%%%%%%%%%%%%%%%%%%%%%%%%%%%%%%%%%%%%%%%%%%%%%%%%%%%%%
\subsection{Construction of the corrector} \label{ssec:constr} 

\vskip.5cm
{\bf Random fields}: we recall that $\Omega=\{0,1\}^{\E_d}$ is the set of sub-graphs of $\Z^d$. 
We shall denote with $\BB$ the set of neighbors of the origin in $\Z^d$. With some abuse of notation, 
we write $\o(b)$ instead of $\o(0,b)$ when $b\in\BB$. 
We use the notation $x.\o$ to denote the natural action of $\Z^d$ on $\Omega$ by translations. 
$\Omega$ is equiped with the product sigma field. 

We endow $\Omega\times\BB$ with the measure $M$ defined by 
\beqnn \int u d M = Q[\sum_{b\in\BB} \o(b) u(\o,b) \1_{\#\CC_0(\o)=\infty}]\,. \eeqnn 
Note that if two random fields $u$ and $v$ coincide in  $L^2(\Omega\times\BB,M)$, then, 
$Q$.a.s. on the event $\#\CC_0(\o)=\infty$, $u(\o,b)=v(\o,b)$ for any $b\in\BB$ 
such that $\o(b)=1$. 

Let $u:\Omega\rightarrow \R$. 
$u$ is said to be {\it local} if it only depends on a finite number of coordinates.  
We associate to $u$ its {\it gradient}: 
$\nabla^{(\o)} u: \Omega\times\BB$ defined by 
\beqnn \nabla^{(\o)} u(\o,b)=u(b.\o)-u(\o)\,.\eeqnn 
Let $L^2_{pot}$ be the closure in $L^2(\Omega\times\BB,M)$ of the set 
of gradients of local fields, and  $L^2_{sol}$ be its orthogonal complement  
in  $L^2(\Omega\times\BB,M)$. 

Fields in $L^2_{pot}$ satisfy a {\it co-cycle} relation: on the event $\#\CC_0(\o)=\infty$, for any 
$u\in L^2_{pot}$ and any closed path in $\CC(\o)$ of the form 
$\gamma=(x_0,x_1,...,x_k)$ with $x_i\sim x_{i+1}$, $\o(x_i,x_{i+1})=1$ and 
$x_0=x_k=0$ then 
$\sum_{i=1}^k u(x_{i-1}.\o,x_i-x_{i-1})=0$.

Let us write down explicitely what it means for a square integrable field $v$ to 
be in $L^2_{sol}$: let $u$ be a local function on $\Omega$. Then 
\beqnn 
&&Q[(\sum_{b\in\BB}\o(b)\, v(\o,b)\,\nabla^{(\o)}\, u(\o,b)\,\1_{\#\CC_0(\o)=\infty}]  
= Q[\sum_{b\in\BB} v(\o,b)\,\nabla^{(\o)} \,u(\o,b) \,\1_{0\in\CC(\o),b\in\CC(\o)}]\\  
&=& Q[\sum_{b\in\BB} v(\o,b)\,( u(b.\o)-u(\o) ) \,\1_{0\in\CC(\o),b\in\CC(\o)}]\,.
\eeqnn 
Using the translation invariance of $Q$, we then get that 
\beqnn 
&&Q[\sum_{b\in\BB} v(\o,b)\, u(b.\o)\, \1_{0\in\CC(\o),b\in\CC(\o)}] 
= Q[\sum_{b\in\BB} v((-b).b.\o,b)\, u(b.\o)\, \1_{0\in\CC(b.\o),-b\in\CC(b.\o)}]\\ 
&=& Q[\sum_{b\in\BB} v((-b).\o,b)\, u(\o)\, \1_{0\in\CC(\o),-b\in\CC(\o)}]  
= Q[\sum_{b\in\BB} v(b.\o,-b)\, u(\o)\, \1_{0\in\CC(\o),b\in\CC(\o)}]\\ 
&=&  Q[\sum_{b\in\BB} \o(b)\, v(b.\o,-b)\, u(\o)\, \1_{\#\CC_0(\o)=\infty}] \,.
\eeqnn 
So that 
\beqnn 
 Q[\sum_{b\in\BB} v(\o,b)\,( u(b.\o)-u(\o) ) \,\1_{0\in\CC(\o),b\in\CC(\o)}]  
= Q[\sum_{b\in\BB} \o(b)\, u(\o)\, (v(b.\o,-b)-v(\o,b))\, \1_{\#\CC_0(\o)=\infty}] \,.
\eeqnn 

Thus we have proved the following integration by parts formula:  
\beq \label{eq:divergence} 
\int v\, \nabla^{(\o)} u\,  dM 
=-Q[n^\o(0)\, u\, \nabla^{(\o)*} v\, \1_{\#\CC_0(\o)=\infty}]\,,
\eeq 
where 
\beqnn 
\nabla^{(\o)*} v(\o)=\frac 1{n^\o(0)} \sum_{b\in\BB} \o(b)\, (v(\o,b)-v(b.\o,-b))\,. 
\eeqnn 
(\ref{eq:divergence}) holds for a square integrable random field  $v$ and any local function $u$. 

As a consequence, taking $v$ to be constant, note that $\int \nabla^{(\o)} u dM =0$ for any local $u$. 
By extension, we will also have $\int  u dM =0$ for any $u\in L^2_{pot}$. 

A square integrable random field $v$ is in $L^2_{sol}$ if 
it satisfies $\nabla^{(\o)*}v=0$ $Q$.a.s. on the set $\#\CC_0(\o)=\infty$. 

\vskip.5cm 

{\bf Definition of the corrector}: let $b\in\BB$. 
Define the random field ${\hat b}(\o,e)=\1_{e=b}-\1_{e=-b}$. 
Let $G_b$ be the unique solution in $L^2_{pot}$ satisfying the equation 
\beq \label{eq:equation} 
{\hat b} +G_b(\o,e)\in L^2_{sol}\,. 
\eeq 
($G_b$ is simply the projection of $-\hat b$ on $L^2_{pot}$.) 

We define the corrector $\chi:\Omega\times\CC_0(\o)\rightarrow \R^d$ by the equation 
\beq \label{eq:equationcorr} 
\chi(\o,x+e)\cdot b-\chi(\o,x)\cdot b=G_b(x.\o,e)\,, 
\eeq 
for any $x\in\Z^d$, $b,e\in\BB$. 
(In this equation $\chi(.)\cdot b$ stands for the usual scalar product of the two $\R^d$ vectors 
$\chi(.)$ and $b$.  
Note that there is no ambiguity because $G_b=-G_{-b}$ as can be directly seen from equation 
(\ref{eq:equation}).) Observe that, unlike $G_b$, the corrector $\chi$ is not an homogeneous field.

The solution to (\ref{eq:equation}) being unique in $L^2_{pot}$, the value 
of $G_b(\o,e)$ is uniquely determined whenever 
$\#\CC_0(\o)=\infty$ and $e\in\BB$ satisfies $\o(e)=1$. 
Therefore $G_b(x.\o,e)$ is well defined $Q$.a.s. on the set $\#\CC_0(\o)=\infty$ 
for any $x$ and $e$ s.t. $x$ and $x+e$ belong to $\CC(\o)$. Thus, if $x$ belongs to $\CC_0(\o)$, 
then the value of 
$\chi(\o,x)-\chi(\o,0)$ can be computed integrating (\ref{eq:equationcorr}) along a path 
in $\CC_0(\o)$ 
from the origin to $x$. 
That this value does not depend on the choice of the path is an immediate consequence 
of the co-cycle relation satisfied by $G_b$. We conclude that $\chi(\o,x)$ is uniquely determined 
by equation (\ref{eq:equationcorr}) up to an additive constant (that might depend on $\o$).

\vskip.5cm 
{\bf The martingale property}: we claim that the random process 
$X(t)+\chi(\o,X(t))$ is a martingale under $P^\o_0$ for $Q$ almost all $\o$ 
s.t. $\#\CC_0(\o)=\infty$. Note that since 
the process $X(t)$, starting from the origin,  never leaves $\CC(\o)$, 
$\chi(\o,X(t))$ is indeed well defined.

We choose $\o$ s.t. $\#\CC_0(\o)=\infty$. 

Since $G_b\in L^2_{pot}$, the co-cycle relation implies that $G_b(\o,e)+G_b(e.\o,-e)=0$ 
for any $e\in\BB$ s.t. $\o(e)=1$. 
Comparing the expression of $\LL^\o$ with the definition of $\nabla^{(\o)*}$, we then see 
that $\LL^\o \chi(\o,x)\cdot b=\frac 1 2\nabla^{(\o)*} G_b(x.\o)$ for any $x\in\CC(\o)$. 

Let $\phi(x)=x+\chi(\o,x)$. 
Noting that ${\hat b}(\o,e)=e\cdot b$ and that 
${\hat b}(e.\o,-e)=-b\cdot e$, we see that 
$\nabla^{(\o)*}{\hat b}(\o)=\frac 2{n^\o(0)}\sum_{e\in\BB} \o(e) b\cdot e$. Therefore 
\beqnn 
&&\LL^\o\phi(\o,x)\cdot b
=\frac 1{n^\o(x)}\sum_{e\in\BB} \o(x,x+e) e\cdot b +\LL^\o\chi(\o,x) \\ 
&=& \frac 1 2\nabla^{(\o)*}{\hat b}(x.\o)+\frac 1 2\nabla^{(\o)*} G_b(x.\o)=0\,. 
\eeqnn 
This last equality holds for any $x\in\CC(\o)$. 
We have proved the martingale property. 

%%%%%%%%%%%%%%%%%%%%%%%%%%%%%%%%%%%%%%%%%%%%%%%%%%%%%%%%%%%%%%%%%%%%%%%%%%%%%%%%%%%%%%%%%%%%%%%%%%%%%%%%%
\vskip.5cm 
{\bf The invariance principle}: 
to each pair of neighbouring points $x,y\in\CC_0(\o)$ such that $\o(x,y)=1$ 
attach a Poisson process of rate $1/n^\o(x)$, say $N^{x,y}_t$, all of them  
independent. 
Let $Y$ be the solution of the equation $Y(0)=0$, 
\beqnn 
dY(t)=\sum_{y\sim Y(t-)} \o(Y(t-),y) (y-Y(t-)) \,dN^{Y(t-),y}_t\,. 
\eeqnn 
Then the law of the random process $(Y(t),t\ge 0)$ is $P^\o_0$. 

Let $\o$ be such that $\#\CC_0(\o)=\infty$. Let $N(t)=Y(t)+\chi(\o,Y(t))$. 
From the previous paragraph, we already know that $N$ is a 
martingale. Its bracket can be computed using It\^o's 
formula. We fix a direction $b\in\BB$. Then: 
\beqnn 
&&d<N\cdot b>(t)\\ 
&=&\frac 1 {n^\o(Y(t-))} \sum_{y\sim Y(t-)} 
     \o(Y(t-),Y(t-)+e)\,(y\cdot b+\chi(\o,y)\cdot b-Y(t-)\cdot b-\chi(\o,Y(t-))\cdot b)^2\, dt\\ 
&=& \frac 1 {n^{Y(t-).\o}(0)} \sum_{e\in\BB} Y(t-).\o(e)\, (e\cdot b+G_b(Y(t-).\o,e))^2 \,dt \,. 
\eeqnn 

Back to the process $X$ and denoting $M(t)=X(t)+\chi(\o,X(t))$, we can equivalently 
write that 
\beqnn 
d<M\cdot b>(t)=\frac 1 {n^{X(t-).\o}(0)} \sum_{e\in\BB} X(t-).\o(e)\, (e\cdot b+G_b(X(t-).\o,e))^2\, dt \,. 
\eeqnn  

Let ${\tilde Q}_0$ be the probability measure 
\beqnn 
{\tilde Q}_0(A)=\frac{\int_A n^\o(0)dQ_0(\o)}{\int n^\o(0)dQ_0(\o)}\,. 
\eeqnn 
$X(t-).\o$ is the process of the environment viewed from the 
particule. The measure ${\tilde Q}_0$ is reversible, invariant and ergodic with respect to 
$X(t-).\o$, see Lemma 4.9 in \cite{kn:DFGW}. 
As a consequence, we get that, $Q$.a.s. on the set $\#\CC_0(\o)=\infty$, 
\beqnn 
\frac {<M\cdot b>(t)} t \rightarrow {\tilde Q}_0( \frac 1 {n^\o(0)} \sum_{e\in\BB} \o(e)\, (e\cdot b+G_b(\o,e))^2 )\,. 
\eeqnn 
Let now $M^\eps(t)=\eps M(t/\eps^2)$. We have proved that, for any $t>0$, as $\eps$ tends to $0$ 
\beqnn 
<M^\eps\cdot b>(t) \rightarrow t {\tilde Q}_0( \frac 1 {n^\o(0)} \sum_{e\in\BB} \o(e)\, (e\cdot b+G_b(\o,e))^2 )\,. 
\eeqnn 

For any  function $f$ that vanishes on the diagonal, the process 
\beqnn 
\sum_{0\leq s\leq t} f(X(s),X(s-)) 
-\int_0^t ds \frac 1 {n^{X(s-).\o}(0)} \sum_{e\in\BB} X(s-).\o(e)\, f(X(s-)+e,X(s-)) 
\eeqnn  
is a martingale. Applying this to 
$f(x,y)=( b\cdot (x+\chi(\o,x))-b\cdot (y+\chi(\o,y)))^2\,   
\1_{\vert b\cdot (x+\chi(\o,x))-b\cdot (y+\chi(\o,y))\vert \geq \eta} $ for some direction $b$ 
and some $\eta>0$, we get that 
\beqnn 
\sum_{0\leq s\leq t} &&( M(s)\cdot b-M(s-)\cdot b)^2 \,\1_{\vert M(s)\cdot b-M(s-)\cdot b\vert\geq\eta} \\ 
&-&\int_0^t ds \frac 1 {n^{X(s-).\o}(0)} \sum_{e\in\BB} X(s-).\o(e)\, ( e\cdot b +G_b(X(s-).\o,e))^2\,  
                     \1_{\vert e\cdot b +G_b(X(s-).\o,e)\vert\geq\eta} 
\eeqnn   
is a martingale. Taking expectations and using the ergodic theorem for the process $X(s-).\o$ 
we get that, on the set $\#\CC_0(\o)=\infty$, 
\beqnn 
&&E^\o_0[\frac 1 t \sum_{0\leq s\leq t} ( M(s)\cdot b-M(s-)\cdot b)^2\, \1_{\vert M(s)\cdot b-M(s-)\cdot b\vert\geq\eta}] \\ 
&=& \frac 1 t  
\int_0^t ds  E^\o_0[ \frac 1 {n^{X(s-).\o}(0)} \sum_{e\in\BB} X(s-).\o(e) ( e\cdot b +G_b(X(s-).\o,e))^2\,  
                     \,\1_{\vert e\cdot b +G_b(X(s-).\o,e)\vert\geq\eta} ]\\ 
&\rightarrow&  
{\tilde Q}_0( \frac 1 {n^{\o}(0)} \sum_{e\in\BB} \o(e)\, ( e\cdot b +G_b(\o,e))^2\,  
                     \1_{\vert e\cdot b +G_b(\o,e)\vert\geq\eta}  )<\infty 
                     \,. \eeqnn     
Then, for any $t>0$ 
\beqnn 
&&E^\o_0[ \sum_{0\leq s\leq t} ( M^\eps(s)\cdot b-M^\eps(s-)\cdot b)^2\, \1_{\vert M^\eps(s)\cdot b-M^\eps(s-)\cdot b\vert\geq\eta}] \\ 
&=& \eps^2\, E^\o_0[ \sum_{0\leq s\leq t/\eps^2} ( M(s)\cdot b-M(s-)\cdot b)^2\, \1_{\vert M(s)\cdot b-M(s-)\cdot b\vert\geq\eta/\eps}] \\ 
&\rightarrow&0\,. 
\eeqnn

From the martingale convergence theorem, Theorem 5.1 part a in \cite{kn:He}, we then deduce that, 
$Q$.a.s. on the set $\#\CC_0(\o)=\infty$, the law of the process 
$\eps X(./\eps^2)+\eps\chi(\o,X(./\eps^2))$ under $P^\o_0$ converges to the law of a 
Brownian motion with a deterministic covariance matrix $A$. 

That $A$ is diagonal is proved in \cite{kn:DFGW}, Theorem 4.7, 3. 
One can argue that $A$ is positive as a consequence of the Gaussian 
lower bounds obtained in \cite{kn:Ba}, but the original proof is given 
in \cite{kn:GM}.  

We therefore conclude that Theorem \ref{theo} will follow if we can prove that, 
$Q$.a.s. on the set $\#\CC_0(\o)=\infty$, for all $t>0$, 
$\eps \chi(\o,X(t/\eps^2))$ converges to $0$ in $P^\o_0$ probability.

\vskip1cm

%%%%%%%%%%%%%%%%%%%%%%%%%%%%%%%%%%%%%%%%%%%%%%%%%%%%%%%%%%%%%%%%%%%%%%%%%%%%%%%%%%%%%%%%%%%%%%%%%%%%%%% 
\subsection{Convergence of the corrector} \label{ssec:conv}

We now check that the contribution of the corrector 
is negligible in the limit i.e. we prove that, for all $t$, 
$\eps\chi(\o,X(t/\eps^2))$ converges to $0$ in $P_0^\o$ probability, 
$Q$.a.s. on the set $\#\CC_0(\o)=\infty$. 
In view of (\ref{eq:Barlow}), it is sufficient to show that 
\beqnn 
\lim_{\eps\rightarrow 0} \eps^d \sum_{y\in\CC_0(\o)\,;\, \vert y\vert\le 1/\eps} 
 \vert \eps\chi(\o,y)\vert^2
=0\, \hbox{\ $Q_0$.a.s.} \,. 
\eeqnn  

Below, we use the Poincar\'e inequality to prove that there exist some constants 
$a_\eps(\o)$ such that 
\beq\label{eq:correct}
\lim_{\eps\rightarrow 0} \eps^d \sum_{y\in\CC_0(\o)\,;\, \vert y\vert\le 1/\eps} 
 \vert \eps\chi(\o,y)-a_\eps\vert^2
=0\, \hbox{\ $Q_0$.a.s.} \,. 
\eeq

As a consequence of (\ref{eq:Barlow}), (\ref{eq:correct}) implies that  
\beqnn 
\lim_{t\rightarrow 0}\lim_{\eps\rightarrow 0} P^\o_0[\vert \eps \chi(\o,X(\frac t{\eps^2}))-a_\eps\vert\ge K]
=0\, \hbox{\ $Q_0$.a.s., and for any $K>0$} \,.    
\eeqnn 

But  the invariance principle for the process 
$\eps X(\frac t{\eps^2})+\eps \chi(\o,X(\frac t{\eps^2}))$ implies that 
\beqnn 
\lim_{t\rightarrow 0}\lim_{\eps\rightarrow 0} P^\o_0[\vert \eps X(\frac t{\eps^2})+\eps \chi(\o,X(\frac t{\eps^2}))\vert\ge K]
=0\, \hbox{\ $Q_0$.a.s., and for any $K>0$} \,,      
\eeqnn 

and (\ref{eq:Barlow}) implies that 
\beqnn 
\lim_{t\rightarrow 0}\lim_{\eps\rightarrow 0} P^\o_0[\vert \eps X(\frac t{\eps^2})\vert\ge K]
=0\, \hbox{\ $Q_0$.a.s., and for any $K>0$} \,.     
\eeqnn 

Therefore 
\beqnn 
\lim_{t\rightarrow 0}\lim_{\eps\rightarrow 0} P^\o_0[\vert \eps \chi(\o,X(\frac t{\eps^2}))\vert\ge K]
=0\, \hbox{\ $Q_0$.a.s., and for any $K>0$} \,.      
\eeqnn 
Thus we see that  $a_\eps$ tends to $0$ and 
\beqnn 
\lim_{\eps\rightarrow 0} \eps^d \sum_{y\in\CC_0(\o)\,;\, \vert y\vert\le 1/\eps} 
 \vert \eps\chi(\o,y)\vert^2
=0\, \hbox{\ $Q_0$.a.s.} \,. 
\eeqnn  

It remains to justify (\ref{eq:correct}).

%%%%%%%%%%%%%%%%%%%%%%%%%%%%%%%%%%%%%%%%%%%%%%%%%%%%%%%%%%%%%%%%%%%%%%%%%%%%%%%%%%%%%%%%%%%%%%%%%%%%%%%%%%%
\vskip.5cm 
{\bf Poincar\'e inequalities}: since $G_b$ is square integrable,  
the spatial ergodic theorem, see \cite{kn:Kr} page 205,  implies that   
$\eps^d \sum_{e\in\BB}\sum_{x\in\CC_0(\o)\,;\, \vert x\vert\le 1/\eps} \o(e) (G_b(x.\o,e))^2$ 
has a $Q$.a.s. finite limit. Therefore 
\beq\label{eq:poinc1}
\limsup_\eps \eps^d \sum_{e\in\BB}\sum_{x\in\CC_0(\o)\,;\, \vert x\vert\le 1/(1-a)\eps} 
x.\o(e)\, (G_b(x.\o,e))^2 <\infty\,,
\eeq  
$Q_0$.a.s. and for any constant $0<a<1$.

We quote from 
\cite{kn:MR}, Theorem 1.3: for  some $\eps>0$  define 
$\CC^\eps$ to be the connected component of the intersection of $\CC_0(\o)$ 
with the box $[-\frac 1 \eps,\frac 1 \eps]^d$ that contains the origin. 
There exists a constant $\beta$ such that 
$Q_0$.a.s. for small enough $\eps$, for any function $u:\CC^\eps\rightarrow\R$ one has 
\beqnn 
\frac 1{\#\CC^\eps} \sum_{x,y\in\CC^\eps} (u(x)-u(y))^2
\le \beta \eps^{-2}\,   
\sum_{x\sim y\in\CC^\eps} \o(x,y)\,(u(x)-u(y))^2\,. 
\eeqnn 

Since $\#\CC^\eps$ is of order $\eps^{-d}$ for small enough $\eps$ and since 
$\CC_0(\o)\cap [-\frac 1 \eps,\frac 1 \eps]^d\subset \CC^{(1-a)\eps}$ for some constant $a$, 
we therefore have a constant $\beta$ such that, 
$Q_0$.a.s. for small enough $\eps$, for any function $u:\CC(\o)\rightarrow\R$ 
\beqnn
\eps^d \,\sum_{x,y\in\CC(\o)\,;\, \vert x\vert\,,\,\vert y\vert \le 1/\eps} (u(x)-u(y))^2
\le \beta \eps^{-2}\,
\sum_{x\sim y\in\CC^{(1-a)\eps} } \o(x,y)\,(u(x)-u(y))^2\,. 
\eeqnn  

We use this last inequality for the function $u(x)=\chi(\o,x)$ to get that, 
\beqnn 
 \eps^d\, \sum_{x,y\in\CC(\o)\,;\, \vert x\vert\,,\,\vert y\vert \le 1/\eps} (\chi(\o,x)-\chi(\o,y))^2
\le \beta \eps^{-2}\, 
\sum_{b\in\BB} \sum_{e\in\BB} \sum_{x\in\CC_0(\o)\,;\, \vert x\vert\le 1/(1-a)\eps} 
x.\o(e)\, (G_b(x.\o,e))^2 \,. 
\eeqnn 

By (\ref{eq:poinc1}) we therefore get that   
\beqnn 
\limsup_\eps { \eps^{2d}\,\sum}_{x,y\in\CC(\o)\,;\, \vert x\vert\,,\,\vert y\vert \le 1/\eps} 
       (\eps\chi(\o,x)-\eps\chi(\o,y))^2<\infty\,, 
       \eeqnn 
$Q_0$.a.s.,  
and thus 
\beq\label{eq:poinc2}  
\limsup_\eps {\eps^d\, \sum}_{x\in\CC(\o)\,;\, \vert x\vert \le 1/\eps} 
       (\eps\chi(\o,x)-a_\eps)^2<\infty\,, 
       \eeq  
$Q_0$.a.s. 
where $a_\eps$ is the mean value of $\eps\chi(\o,x)$ on the set 
$\{x\in\CC(\o)\,;\, \vert x\vert \le 1/\eps\}$.

%%%%%%%%%%%%%%%%%%%%%%%%%%%%%%%%%%%%%%%%%%%%%%%%%%%%%%%%%%%%%%%%%%%%%%%%%%%%%%%%%%%%%%%%%%%%%%%%%%%%%%%%%%% 

\def\bmu{{\cal P}}

\vskip.5cm 
{\bf Two scale convergence}: we first introduce some notation. 
Let $G=]-1,1[^d$. 
For $\o\in\Omega$ and $\eps>0$, we define the measures  
\beqnn 
\mu_\o=\sum_{z\in \CC(\o)} n^\o(z)\, \delta_z \, ,\, 
\mu^\eps_\o=\eps^d \sum_{z\in \CC(\o)} n^\o(z)\, \delta_{\eps z}\,.
\eeqnn  
Given a direction $e\in\BB$, the gradient of a function $\phi:\R^d\rightarrow\R$ is  
\beqnn 
\nabla^\eps_e\phi(z)=\frac 1\eps (\phi(z+\eps e)-\phi(z))\,. 
\eeqnn 
Let us now choose $b_0\in\BB$ and let 
\beqnn 
\psi^\eps(\o,z)=(\eps\chi(\o,\frac 1\eps z)-a_\eps)\cdot b_0\,.
\eeqnn 
Thus $\psi^\eps$ is well defined for $z\in\eps\,\CC_0(\o)$. 
From the definition of $\chi$, we have 
\beqnn 
\nabla^\eps_e\psi^\eps(\o,z)=G_{b_0}(\frac 1\eps z.\o,e)\,,
\eeqnn 
for $z\in \eps\,\CC_0(\o)$. 

{\it Keep in mind that for $z'\in\Z^d$, the expression  $z'.\o$ denotes the graph 
obtained by translating $\o$ by $z'$. In 
particular, for $z\in\eps\Z^d$, then $\frac 1\eps z.\o(e)$ is either $0$ or $1$, depending on wether the edge 
$(z,z+e)$ belongs to $\o$ or not. We sometimes prefer the notation $(\frac z\eps).\o(e)$ in order to avoid 
possible confusion.} 

In our new notation, (\ref{eq:poinc1}) and (\ref{eq:poinc2}) now read: 
\beq \label{eq:poinc11} 
C_1(\o)=\sup_{e\in\BB}\sup_\eps \int_G (\frac z\eps).\o(e)\,(\nabla^\eps_e\psi^\eps(\o,z))^2\, d\mu^\eps_\o(z)<\infty\,,
\eeq  
and 
\beq \label{eq:poinc22}  
C_2(\o)=\sup_\eps \int_G (\psi^\eps(\o,z))^2\, d\mu^\eps_\o(z)<\infty\,, 
\eeq  
for $Q_0$ almost any $\o$. 
For further reference, let us call $\Omega_1$ the set of $\o$'s such that $0\in\CC(\o)$, 
$C_1(\o)<\infty$ and $C_2(\o)<\infty$ and observe that $Q_0(\Omega_1)=1$. 

Define the measure 
\beqnn 
\bmu(A)=Q[ \1_A(\o)\, n^\o(0)\, \1_{0\in\CC(\o)}]. 
\eeqnn 
According to the ergodic theorem, for any smooth function $\phi\in C^\infty(G)$ 
and any $u\in L^1(\O,\bmu)$ we have 
\beq \label{eq:ergo} 
\int_G \phi(z)\, u(\frac 1\eps z.\o)\, d\mu^\eps_\o(z) 
\rightarrow (\int_G\phi(z)\,dz) (\int_\Omega u(\o')\,d\bmu(\o'))\,,
\eeq 
$Q_0$.a.s. 

We endow $\Omega$ with its natural (product) topology to turn it into a compact 
space. We will use the notation $C(\O)$ for continuous real valued functions defined on $\Omega$. 
Using standart separability arguments, we see that (\ref{eq:ergo}) holds simultaneously 
for any $\phi\in C^\infty(G)$ and  $u\in  C(\O)$ on a set of full $Q_0$ measure. 
More precisely, let $\O_2$ be the set of $\o$'s such that $0\in \CC(\o)$ and, for any functions 
$\phi\in C^\infty(G)$ and  $u\in C(\O)$ one has: 
\beq \label{eq:ergo11}  
\int_G \phi(z) u(\frac 1\eps z.\o) d\mu^\eps_\o(z) 
\rightarrow (\int_G\phi(z)dz) (\int_\Omega u(\o')d\bmu(\o'))\,, 
\eeq 
and, for any $e\in\BB$,   
\beq \label{eq:ergo22} 
\int_G \phi(z) u(\frac 1\eps z.\o) G_{b_0}(\frac 1\eps z.\o,e) d\mu^\eps_\o(z) 
\rightarrow (\int_G\phi(z)dz) (\int_\Omega u(\o')G_{b_0}(\o',e)d\bmu(\o'))\,.  
\eeq  
Then $Q_0(\Omega_2)=1$. 
Finally let $\Omega_0=\Omega_1\cap\Omega_2$. 

In the sequel $\a$ will denote an element of $\Omega_0$. 
Consider the family of linear functionals 
\beqnn 
L^{\eps,\a}(u,\phi)=\int_G \phi(z) \psi^\eps(\a,z) u(\frac 1\eps z.\a) d\mu^\eps_\a(z)\,.
\eeqnn 

Using the Cauchy-Schwartz inequality, we get that 
\beqnn 
(L^{\eps,\a}(u,\phi))^2 \leq 
\int_G (\psi^\eps(\a,z))^2\, d\mu^\eps_\a(z)\, \int_G \phi(z)^2\, u(\frac 1\eps z.\a)^2 \,d\mu^\eps_\a(z)\,. 
\eeqnn 
From (\ref{eq:poinc22}) and (\ref{eq:ergo11}) we deduce 
that for $\phi\in C^\infty(G)$ and  $u\in  C(\O)$ 
\beqnn 
\limsup_\eps (L^{\eps,\a}(u,\phi))^2 \leq 
C_2(\a) \int_G \phi(z)^2\, dz \int_\Omega u(\o)^2\, d\bmu(\o)\,.
\eeqnn 
Therefore, 
up to extracting a sub-sequence, we can assume that for any 
smooth $\phi$ and any continuous $u\in C(\O)$, $L^{\eps,\a}(u,\phi)$ has a limit  say 
$L^\a(u,\phi)$ where $L^\a$ is a linear functional satisfying 
\beqnn 
(L^\a(u,\phi))^2\, \leq 
C_2(\a) \int_G \phi(z)^2\,  dz \int_\Omega u(\o)^2\,  d\bmu(\o)\,.
\eeqnn 
Thus $L^\a$ can be extended as a continuous 
linear functional on $L^2(\Omega\times G,d\bmu\times dx)$
and, by Riesz's theorem, there exists a function 
$v^\a\in L^2(\Omega\times G,d\bmu\times dx)$ such that 
\beqnn L^\a(u,\phi)=\int_G \phi(z)\, dz\int_\Omega u(\o)\, v^\a(\o,z)\, d\bmu(\o)\,.
\eeqnn 

Let us summarize the preceeding discussion: we have proved that, up to extracting a sub-sequence, 
for $\phi\in C^\infty(G)$ and $u\in C(\Omega)$, 
\beq \label{eq:2scale} 
\int_G \phi(z)\, \psi^\eps(\a,z)\, u(\frac 1\eps z.\a)\, d\mu^\eps_\a(z) 
\rightarrow \int_G \phi(z)\, dz\int_\Omega u(\o)\, v^\a(\o,z)\, d\bmu(\o)\,. 
\eeq 

We will prove the following 

\begin{lm} \label{lemma} 
For any $\a\in\Omega_0$, $v^\a(\o,z)=0$ for Lebesgue almost any 
$z\in G$ and $\bmu$ almost any $\o$.
\end{lm}  

As a consequence of this Lemma, we have that for $Q_0$ almost any $\a$, for 
any function $\phi\in C^\infty(G)$, 
\beqnn 
\int_G \phi(z)\, \psi^\eps(\a,z)\,  d\mu^\eps_\a(z) 
\rightarrow 0\,.
\eeqnn 

Since we also have uniform bounds on the $L^2$ norm of $\psi^\eps$, see (\ref{eq:poinc22}), we deduce that, 
for any rectangle $A\subset G$, 
\beqnn 
\int_A \psi^\eps(\a,z)\,  d\mu^\eps_\a(z) 
\rightarrow 0\,.
\eeqnn 

We conclude that, for any rectangle $A\subset [-1,1]^d$, $Q_0$.a.s. 
\beq\label{eq:local} 
\eps^d\, \sum_{x\in\CC(\o)\,;\, \eps x\in A}(\eps\chi(\o,x)-a_\eps)  
\rightarrow 0\,.
\eeq

\begin{rmk} The content of this part of the paper, including the proof of the Lemma in the next section, 
should be compared with the results of \cite{kn:JP}. The convergence in (\ref{eq:2scale}) is known as 
'two-scale convergence'. The only difference between our setting and \cite{kn:JP} is the discrete 
nature of the grid whereas continuous diffusions are considered in \cite{kn:JP}. 

It is also possible to directly apply the results of \cite{kn:JP} to justify Lemma (\ref{lemma}). 
We refer the interested reader to the first version of the present paper on the Arxiv for details. 
Here, we prefered to give a more self-contained approach but most of the arguments are mere copies 
of the proofs in  \cite{kn:JP} with some minor simplifications due to the fact that, for instance, 
the Palm measure $\bmu$ is explicit and absolutely continuous w.r.t. $Q$. 
\end{rmk}

%%%%%%%%%%%%%%%%%%%%%%%%%%%%%%%%%%%%%%%%%%%%%%%%%%%%%%%%%%%%%%%%%%%%%%%%%%%%%%%%%%%%%%%%%%%%%%%%%%%%%%%%%%% 

\vskip.5cm 
{\bf Proof of Lemma (\ref{lemma})}: the proof is in three steps. 
Throughout the following proof, $\phi$ is always assumed to be in $C_o^\infty(G)$, the 
space of smooth functions with compact support in $G$.  

{\bf Step 1}: we check the integration by parts formula: 
\beq \label{eq:intpart} 
\int_G \phi(z)\, \nabla^{(\o)*} u(\frac 1 \eps z.\a)\, d\mu^\eps_\a(z) 
=-\eps\int_G \frac 1{n^\a(\frac 1\eps z)} \sum_{e\in\BB} 
u(\frac 1 \eps z.\a,e) \,(\frac z\eps) .\a(e) \, 
\nabla^\eps_e\phi(z)\, d\mu^\eps_\a(z)\,,
\eeq 
where  $u$ is any function defined on $\Omega\times\BB$ and 
$\eps$ is small enough (depending on the support of $\phi$):  

\beqnn 
&& \int_G \phi(z) \nabla^{(\o)*}\, u(\frac 1 \eps z.\a)\, d\mu^\eps_\a(z) 
= \eps^d \sum_{x\in \CC(\a)} \phi(\eps x)\, \nabla^{(\o)*}u(x.\a)\, n^\a(x) \\
&=& \eps^d \sum_{x\in \CC(\a)} \phi(\eps x)\, \sum_{e\in\BB} x.\a(e)\, (u(x.\a,e)-u(x.e.\o,-e)) \,n^\a(x)\,. 
\eeqnn 
But 
\beqnn 
\sum_{x\in \CC(\a)} \phi(\eps x)\, \sum_{e\in\BB} x.\a(e)\, u(x.e.\a,-e)\, n^\a(x) 
= \sum_{x'\in \CC(\a)} \sum_{e'\in\BB} x'.\a(e')\, u(x'.\a,e')\, \phi(\eps x'+\eps e')\, n^\a(x')\,, 
\eeqnn 
with the change of variables $x'=x+e$ and $e'=-e$. 
Putting the last two equalities together, one gets (\ref{eq:intpart}). 
Observe that boundary terms vanish because 
$\phi$ has compact support and $\eps$ is small enough. 

{\bf Step 2}: we prove that $v^\a(\o,z)$ does not depend on $\o$ i.e. that 
$Q_0$.a.s. 
\beq \label{eq:stp2} 
v^\a(\o,z)=\frac {\int v^\a(\o',z)\,d\bmu(\o')}{\int d\bmu(\o')}=v^\a(z)\,.
\eeq 

Indeed, let $u$ be continuous on $\Omega\times\BB$  and $\phi\in C_o^\infty(G)$ 
and use (\ref{eq:2scale}) and the integration by parts 
formula (\ref{eq:intpart}) to get that 
\beqnn 
&& \int_G \phi(z)\, dz\int_\Omega  v^\a(\o,z)\, \nabla^{(\o)*} u(\o)\, d\bmu(\o)\\ 
&=&\lim_\eps \int_G \phi(z)\, \psi^\eps(\a,z)\, \nabla^{(\o)*}u(\frac 1\eps z.\a)\, d\mu^\eps_\a(z)\\ 
&=&\lim_\eps -\eps\int_G \frac 1{n^\a(\frac 1\eps z)} \sum_{e\in\BB}  
u(\frac 1 \eps z.\a,e) \,(\frac z\eps) .\a(e) \, 
\nabla^\eps_e(\psi^\eps(\a,.)\phi)(z)\, d\mu^\eps_\a(z)\,. 
\eeqnn 
Since $u$ is continuous, it is bounded. 
Note that $(\frac z\eps) .\a(e)\leq n^\a(\frac 1\eps z)$. 
Besides, 
\beqnn \limsup_\eps 
\int_G (\frac z\eps) .\a(e) \, (\nabla^\eps_e(\psi^\eps(\a,.)\phi)(z))^2\, d\mu^\eps_\a(z)
\leq C_1(\a) \Vert\phi\Vert_\infty+C_2(\a)\Vert\nabla\phi\Vert_\infty<\infty\,.
\eeqnn 
We conclude that, as $\eps$ tends to $0$, the expression 
\beqnn \int_G \frac 1{n^\a(\frac 1\eps z)} 
u(\frac 1 \eps z.\a,e) \,(\frac z\eps) .\a(e) \, 
\nabla^\eps_e(\psi^\eps(\a,.)\phi)(z)\, d\mu^\eps_\a(z)\eeqnn remains bounded 
and therefore 
\beqnn \lim_\eps -\eps\int_G \frac 1{n^\a(\frac 1\eps z)} \sum_{e\in\BB}  
u(\frac 1 \eps z.\a,e) \,(\frac z\eps) .\a(e) \, 
\nabla^\eps_e(\psi^\eps(\a,.)\phi)(z)\, d\mu^\eps_\a(z)=0\eeqnn and 
\beqnn \int_G \phi(z)\, dz\int_\Omega  v^\a(\o,z)\, \nabla^{(\o)*} u(\o) \,d\bmu(\o)=0\,.\eeqnn

By (\ref{eq:divergence}) we also have 
\beqnn 
&& \int_\Omega  v^\a(\o,z)\, \nabla^{(\o)*} u(\o)\, d\bmu(\o) 
=\int_\Omega  v^\a(\o,z)\, \nabla^{(\o)*} u(\o)\, n^\o(0)\, \1_{0\in\CC(\o)}\, dQ(\o)\\
&=& -\int u \nabla^{(\o)} v^\a(.,z)\, dM \,.  
\eeqnn 
Thus we have proved that 
\beqnn 
\int_G \phi(z)\, dz \int u \,\nabla^{(\o)} v^\a(.,z)\, dM =0\,,
\eeqnn 
for any $\phi\in C_o^\infty(G)$ and continuous $u$. We deduce that $Q_0$.a.s., for any $b\in\BB$ such that 
$\o(b)=1$ and for Lebesgue almost any $z$, then $v^\a(b.\o,z)=v^\a(\o,z)$. 
Integrating this equality on a path between $0$ and $x\in\CC_0(\o)$, we then get that 
$Q_0$.a.s. for any $x\in\CC(\o)$ and for Lebesgue almost any $z$, then 
$v^\a(\o,z)=v^\a(x.\o,z)$. 
Therefore, since $\mu^\eps_\o$ only charges $\CC(\o)$, the ergodic theorem yields: 
\beqnn  
v^\a(\o,z)=\frac{\int_G v^\a(\frac 1\eps z'.\o,z)\, d\mu^\eps_\o(z')}{\int_G d\mu^\eps_\o(z')}  
\rightarrow  \frac {\int v^\a(\o',z)\,d\bmu(\o')}{\int d\bmu(\o')}=v^\a(z)\,,
\eeqnn 
$Q_0$.a.s in $\o$ and for Lebesgue almost any $z\in G$.

{\bf Step 3}: we now prove that $v^\a(z)$ does not depend on $z$. 
To this end, we first prove that for any smooth $\phi\in C_o^\infty(G)$ and any continuous 
$u\in L^2_{sol}$ we have 
\beq \label{eq:stp3} 
\sum_{e\in\BB} (\int_G dz\, v^\a(z) \,\nabla\phi(z) \cdot e) 
(\int_\Omega  {\tilde u}(\o,e) \,d\bmu(\o))=0\,,
\eeq 
where ${\tilde u}(\o,e)=\frac{\o(e)}{n^\o(0)} u(\o,e)$. 

We have: 
\beqnn 
&& (\int_G dz\, v^\a(z) \nabla\phi(z) \cdot e)   
(\int_\Omega  {\tilde u}(\o,e) \,d\bmu(\o))\\  
&= &\lim_\eps \int_G (\nabla\phi(z)\cdot e)\, \psi^\eps(\a,z) \,{\tilde u}(\frac 1\eps z.\a,e)\,d\mu^\eps_\a(z)  
= \lim_\eps \int_G \nabla^\eps_e\phi(z)  \,\psi^\eps(\a,z) \,{\tilde u}(\frac 1\eps z.\a,e)\,d\mu^\eps_\a(z)\\
&=& \lim_\eps \int_G \nabla^\eps_e(\phi  \psi^\eps(\a,.))(z) \,{\tilde u}(\frac 1\eps z.\a,e)\,d\mu^\eps_\a(z)
-\lim_\eps \int_G \phi(z)\,\nabla^\eps_e(\psi^\eps(\a,.))(z) \,{\tilde u}(\frac 1\eps z.\a,e)\,d\mu^\eps_\a(z)
\,,
\eeqnn 
where we used (\ref{eq:2scale}) in the first equality  
and the regularity of $\phi$ for the second and third equalities. 
Using integration by parts and the definition of $\tilde u$, we get   
\beqnn 
&&\sum_{e\in\BB} 
\int_G \nabla^\eps_e(\phi  \psi^\eps(\a,.))(z) \,{\tilde u}(\frac 1\eps z.\a,e)\,d\mu^\eps_\a(z)\\ 
&=&\sum_{e\in\BB} 
\int_G \nabla^\eps_e(\phi  \psi^\eps(\a,.))(z) \,\frac{(\frac z\eps) .\a(e)}{n^\a(\frac 1\eps z)} 
 \,u(\frac 1\eps z.\a,e)\,d\mu^\eps_\a(z)\\ 
&=&-\frac 1\eps \int_G \phi(z) \, \psi^\eps(\a,z)\, 
\nabla^{(\o)*} u(\frac 1\eps z.\a)\,d\mu^\eps_\a(z)=0,
\eeqnn 
since $u\in L^2_{sol}$ and therefore $\nabla^{(\o)*} u=0$. 

We now turn to the second term. Keep in mind that 
$ \nabla^\eps_e(\psi^\eps(\a,.))(z) =G_{b_0}(\frac 1\eps z.\a,e)$. Thus, as an application 
of (\ref{eq:ergo22}), 
\beqnn 
\lim_\eps \int_G \phi(z)\,\nabla^\eps_e(\psi^\eps(\a,.))(z) \,{\tilde u}(\frac 1\eps z.\a,e)\,d\mu^\eps_\a(z)
=(\int_G \phi(z)\,dz) (\int_\Omega G_{b_0}(\o,e)\,{\tilde u}(\o,e) \,d\bmu(\o))\,.
\eeqnn 
Replacing ${\tilde u}$ and $\bmu$ by their definition, we also have
\beqnn 
&&\sum_{e\in\BB} \int_\Omega G_{b_0}(\o,e)\,{\tilde u}(\o,e) \,d\bmu(\o)  
=\sum_{e\in\BB} \int_\Omega G_{b_0}(\o,e)\, \o(e)\, u(\o,e)\, \1_{0\in\CC(\o)}\, dQ(\o)\\ 
&=& \int G_{b_0} u \,dM=0\,, 
\eeqnn 
since $G_{b_0}\in L^2_{pot}$ and $u\in L^2_{sol}$. We conclude that 
(\ref{eq:stp3}) holds.

(\ref{eq:stp3}) was proved for any continuous $u\in L^2_{sol}$. By density it also holds 
for any $u\in L^2_{sol}$.  

It remains to check the following fact: for any direction $e\in\BB$, there exists 
$u\in L^2_{sol}$ such that \\ 
$\int_\Omega  {\tilde u}(\o,e) \,d\bmu(\o)\not=0$. 
Indeed, first note that 
\beqnn 
\int_\Omega  {\tilde u}(\o,e)\,d\bmu(\o)
=\int_{0\in\CC(\o)} \o(e)\, u(\o,e)\, dQ(\o)\,.
\eeqnn 
Define the random field $\tilde e$ by  ${\tilde e}(\o,b)=\1_{b=e}$. ($e$ is kept fixed.) 
Let $G$ be the orthogonal projection of $-\tilde e$ on $L^2_{pot}$ and let  
$u=G+{\tilde e}\in L^2_{sol}$. We write that $u$ and $u-{\tilde e}=G$ are orthogonal: 
\beqnn 
\int_{0\in\CC(\o)} \o(e)\, u(\o,e)\, dQ(\o)  
=\int u{\tilde e}\, dM=\int u^2\, dM\not=0\,,
\eeqnn 
because ${\tilde u}\notin L^2_{sol}$.

%Since $G\in L^2_{pot}$, we have 
%\beqnn 
%\int_\Omega  (\o(e)G(\o,e)+\o(-e)G(\o,-e)) dQ(\o)=0\,.
%\eeqnn 

Thus we can deduce from (\ref{eq:stp3}) that 
$\int_G dz\, v^\a(z) \,\nabla\phi(z) \cdot e =0$ for any smooth $\phi$ and any direction $e$. 
Therefore $v^\a$ is Lebesgue almost surely constant. 

{\bf Conclusion of the proof of Lemma (\ref{lemma})}: since $\psi^\eps$ has vanishing mean 
on $G$ - Remember this is the way we chose $a_\eps$ - then $v^\a$ also has vanishing mean in $G$. 
And since, by steps 2 and 3, $v^\a$ is almost surely constant, we must have that 
$Q_0$.a.s. and for Lebesgue almost any $z$, $v^\a(\o,z)=0$. 
\qed

%%%%%%%%%%%%%%%%%%%%%%%%%%%%%%%%%%%%%%%%%%%%%%%%%%%%%%%%%%%%%%%%%%%%%%%%%%%%%%%%%%%%%%%%%%%%%%%%%%%%%%%%%%%
\vskip.5cm 
{\bf Scaling and strong $L^2$ convergence of $\chi$}: 

To conclude the proof of the Theorem, we still have to prove the strong $L^2$ convergence 
in (\ref{eq:correct}). It will be a consequence of the weak convergence (\ref{eq:local}) and 
of a scaling argument. 

We choose a parameter $\delta>0$. We chop the box $[-1,1]^d$ into smaller 
boxes of side length of order $\delta$: for $z\in\delta \Z^d$ s.t. $\vert z\vert\le 1$, let 
$B_z$ (resp. $C_z$) be the box of center $z$ and  side length $M\delta$ 
(resp. side length $\delta$). $M$ is a constant whose value will be chosen later. 
For $\eps>0$, we use the notation 
$B_z(\eps)=(\frac 1 \eps B_z )\cap \Z^d$ and $C_z(\eps)=(\frac 1 \eps C_z )\cap \Z^d$. 

The following version of the Poincar\'e inequality is proved in 
\cite{kn:Ba}, see Definition 1.7, Theorem 2.18, Lemma 2.13 and Proposition 2.17: 
there exist constants $M>1$ and $\beta$ such that 
$Q_0$.a.s. for any $\delta>0$, for small enough $\eps$, for any  
$z \in\delta \Z^d$ s.t. $\vert z\vert\le 1$ and for 
any function $u:\Z^d\rightarrow\R$ one has 
\beqnn 
\frac 1{\# C_z(\eps)} \sum_{x,y\in\CC(\o)\cap C_z(\eps)} (u(x)-u(y))^2
\le \beta \delta^2 \eps^{-2}  
\sum_{x\sim y\in\CC(\o)\cap B_z(\eps)} \o(x,y)\,(u(x)-u(y))^2\,. 
\eeqnn 
 
We use this inequality for the function $\eps\chi$, to get that 
\beqnn 
\frac 1{\# C_z(\eps)} \sum_{x,y\in\CC(\o)\cap C_z(\eps)} (\eps\chi(\o,x)-\eps\chi(\o,y))^2
\le \beta \delta^2   
\sum_{x\in\CC(\o)\cap B_z(\eps)} \sum_{b\in\BB} \sum_{e\in\BB}\o(x,y)\,(G_b(\o,e))^2\,. 
\eeqnn    

Denoting with $a_\eps(z)$ the mean value of $\eps\chi(\o,.)$ on the 
set $\CC(\o)\cap C_z(\eps)$, we get that for all $z$, 
\beqnn 
\sum_{x\in\CC(\o)\cap C_z(\eps)} (\eps\chi(\o,x)-a_\eps(z))^2
\le \beta \delta^2   
\sum_{x\in\CC(\o)\cap B_z(\eps)} \sum_{b\in\BB} \sum_{e\in\BB}\o(x,y)\,(G_b(\o,e))^2\,,    
\eeqnn  
and summing over all values of $z$, 
\beqnn 
\sum_z \sum_{x\in\CC(\o)\cap C_z(\eps)} (\eps\chi(\o,x)-a_\eps(z))^2
\le \beta \delta^2   
\sum_{x\in\CC(\o)\,;\, \vert x\vert \le 1/\eps} \sum_{b\in\BB} \sum_{e\in\BB}\o(x,y)\,(G_b(\o,e))^2\,.   
\eeqnn  

(Remember that the value of $\beta$ is allowed to change from line to line.) 
Multiplying by $\eps^d$ and  applying the spatial ergodic theorem as before,  we get that  
\beqnn 
\limsup_\eps 
\sum_z \eps^d  \sum_{x\in\CC(\o)\cap C_z(\eps)} (\eps\chi(\o,x)-a_\eps(z))^2
\le \beta \delta^2 \sum_{b\in\BB}\int (G_b)^2\, dM\,.    
       \eeqnn   

On the other hand, it follows from (\ref{eq:local}) that, for any $z$, 
$a_\eps(z)-a_\eps$ converges to $0$. Therefore we must also have 
\beqnn 
\limsup_\eps \eps^d 
\sum_z \sum_{x\in\CC(\o)\cap C_z(\eps)} (\eps\chi(\o,x)-a_\eps)^2
\le \beta \delta^2 \sum_{b\in\BB}\int (G_b)^2\, dM\,,    
       \eeqnn 
and 

 \beqnn 
\limsup_\eps {\eps^d \sum}_{x\in\CC(\o)\,;\, \vert x\vert \le 1/\eps} 
       (\eps\chi(\o,x)-a_\eps)^2\le \beta\delta^2 \sum_{b\in\BB}\int (G_b)^2\, dM\,,   
       \eeqnn 

and, since this holds for any $\delta>0$, we deduce that 
\beqnn 
{\eps^d \sum}_{x\in\CC(\o)\,;\, \vert x\vert \le 1/\eps} 
       (\eps\chi(\o,x)-a_\eps)^2\rightarrow 0\,, 
       \eeqnn $Q_0$.a.s. 
\qed

\vfill\eject


\begin{thebibliography}{xxxxxx 89} 

\bibitem{kn:Ba} Barlow, M.T.~(2004)\\ 
                Random walks on supercritical percolation clusters\\  
                {\em Ann.~Probab.}~{\bf 32}, 3024-3084. 
                
\bibitem{kn:BB}   Berger, N., Biskup, M.~(2005)\\  
		Quenched invariance principle for simple random walk 
		on percolation clusters. Preprint.\\              
                http://front.math.ucdavis.edu/math.PR/0503576. 

\bibitem{kn:DFGW} De Masi, A.,  Ferrari, P., Goldstein, S.,  Wick, W.D.~(1989)\\
                  An invariance principle for reversible Markov 
                  processes. Applications to random motions in random environments\\
                 {\em Journ.~Stat.~Phys.}~{\bf 55}~(3/4),  787-855.  
                 
\bibitem{kn:EK} Ethier, S.N., Kurtz, T.G.~(1986)\\ 
		{\em Markov processes}\\ 
		John Wiley, New York.  

		
\bibitem{kn:G}  Grimmett, G.~(1999)\\
                {\em Percolation}\\
                Springer-Verlag, Berlin (Second edition).  
                
\bibitem{kn:GM} Grimmett, G., Marstrand, J.~(1990)\\ 
		The supercritical phase of percolation is well behaved\\                  
		{\em Proc.~Royal~Society~(London)} Ser. A. {\bf 4306}, 429-457. 

\bibitem{kn:He}  Helland, I~(1982)\\ 
		Central limit theorems for martingales with dicrete or continuous time\\ 
		{\em Scand. Journ. Stat.}~{\bf 9}, 79-94.   
		
\bibitem{kn:JKO} Jikov, V.V., Kozlov, S.M., Oleinik, O.A.~(1994)\\
		Homogenization of differential operators and integral functionals.\\ 
		 Springer-Verlag, Berlin. 		
                
\bibitem{kn:JP}  Jikov, V.V., Piatnitski, A.L.~(2005)\\                 
                Homogenization of random singular structures and measures.\\ 
                Preprint (in Russian). 
                
\bibitem{kn:Ko} Kozlov, S.M.~(1985)\\ 
                The method of averaging and walks in inhomogeneous environments\\
                {\em Russian~Math.~Surveys}~{\bf 40}~(2), 73-145.      

    
                               
\bibitem{kn:Kr} Krengel, U~(1985)\\ 
		{\em Ergodic theorems}\\ 
		Walter de Gruyter, Berlin.  

\bibitem{kn:MR} Mathieu, P., Remy, E.~(2004)\\
                Isoperimetry and heat kernel decay on percolations clusters\\ 
                {\em Ann.~Probab.}~{\bf 32}, 100-128. 
                
\bibitem{kn:SS} Sidoravicius, V., Sznitman, A-S.~(2004)\\ 
                Quenched invariance principles for walks on clusters of percolation 
                or among random conductances\\
                {\em Prob.~Th.~Rel.~Fields}~{\bf 129}, 219-244.     

    






\end{thebibliography}
\end{document}